\documentclass[12pt,reqno]{article}

\usepackage[usenames]{color}
\usepackage{amssymb}
\usepackage{graphicx}
\usepackage{amscd}

\usepackage[colorlinks=true,
linkcolor=webgreen,
filecolor=webbrown,
citecolor=webgreen]{hyperref}

\definecolor{webgreen}{rgb}{0,.5,0}
\definecolor{webbrown}{rgb}{.6,0,0}

\usepackage{color}
\usepackage{fullpage}
\usepackage{float}

\usepackage{graphics,amsmath,amssymb}
\usepackage{amsthm}
\usepackage{amsfonts}
\usepackage{latexsym}
\usepackage{epsf}
\usepackage{verbatim}

\usepackage{tikz}
\usepackage{forest}

\usepackage{mathtools}

\usepackage{soul}

\newcounter{puzzle}
\newenvironment{puzzle}[1][]{\refstepcounter{puzzle}
\begin{quote} \textbf{Puzzle~\thepuzzle. #1}}{\end{quote}}

\usepackage[title]{appendix}

\theoremstyle{plain}

\theoremstyle{definition}

\theoremstyle{remark}

\begin{document}

\title{It's Common Knowledge}
\author{Matvey Borodin \and Aidan Duncan \and Joshua Guo \and Kunal Kapoor \and Tanya Khovanova \and Anuj Sakarda \and Jerry Tan \and Armaan Tipirneni \and Max Xu \and Kevin Zhao}
\date{}

\maketitle

\begin{abstract}
We discuss some old common knowledge puzzles and introduce a lot of new common knowledge puzzles.
\end{abstract}


\section{Do you know your hat color?}\label{sec:intro}

Here is a very old puzzle about hats:

\begin{puzzle}[The intro puzzle. The Evidently Evil Emperor and the Three Shrewd Sages.]\label{puzzle:intro}
An evil emperor puts his three sages Alice, Bob, and Charlie to a test. He shows them six hats, three blue hats and three red hats. The sages close their eyes, and the emperor puts a hat on each of their heads. They open their eyes and see the hats on the other sages but not on themselves. The emperor says that at least one of them has a red hat. If at least one of them can figure out the color of their hat, they all can go free, otherwise all their heads are chopped off. What will happen?
\end{puzzle}

Here is the solution. If one sage, let's say Alice, sees two blue hats she knows that her hat is red and announces it. If no one announces the color, that means there is more than one red hat. After some time, if Alice sees exactly one red hat, and neither Bob nor Charlie announce their hat color she realizes that she also has a red hat. After some more time, if no one announces the color, there must be three red hats.

We wanted to see how this works in practice and recreated the set up of the puzzle as an experiment on our friends Alex, Betty and Claire. We put the only blue hat on Alex, the fastest thinker in the group. Not only is Alex a fast thinker, but he also thinks that Betty and Claire are as fast as he is. While Betty was distracted by a cat and Charlie was scratching his head, Alex decided that enough time has passed and said that his hat is red. Oops.

We see that time is an important part of this puzzle. In the solution to the puzzle above, it is not clear how much time the sages should wait for the other sages' announcements. If all of them think at the same speed, the strategy works. However, in reality, they should decide on a fixed amount of time, say one minute, that they will wait for each round of announcements.

We did another experiment with members of our group. Tanya put red and blue hats on all nine students, and said that there is at least one student with a red hat. Then Jerry clapped to start the game and continued clapping every minute. After six claps, all six people with red hats announced that they had a red hat.

Here is their reasoning. After the first clap, if someone sees all blue hats, he can figure out that his has a red hat. After the second clap, everyone knows that there is more than one red hat. Thus, a person who sees exactly one red hat knows that he must have a red hat. If nothing happens before the third clap, everyone knows that there are at least three red hats. Continuing in a similar manner, as there were six hats, after six claps, all the bearers of red hats know that they had red hats. After that, everyone else will know that they have blue hats.

The most interesting part of this puzzle is that Tanya's announcement at the beginning that there was at least one red hat did not tell any of the students anything that they did not already know because every student saw either five or six red hats. So Tanya didn't say anything new. Or did she? Without her announcement, the students wouldn't have been able to figure out their hat colors. How did Tanya's statement help?

Here is the reasoning. 

Suppose only Matvey has a red hat, then he doesn't know that there is at least one red hat, so Tanya's information is new to him and allows him to figure out his hat color. Now suppose only two students, Anuj and Kunal, have red hats. Then everyone knows that there is at least one red hat. But Anuj doesn't know that Kunal knows that there is at least one red hat. If Anuj's hat were blue, then Kunal wouldn't know that there is at least one red hat without Tanya's announcement. After the announcement, since Kunal doesn't say that he knows his hat color after the first clap, Anuj must also have a red hat. If there are three students with red hats, then Tanya's announcement tells the people with red hats that everyone knows that everyone knows that everyone knows that there is at least one red hat.

A fact is called \textit{mutual knowledge} if all the players know the fact. In our experiment with six red hats, the fact that there were at least five red hats on people's heads was mutual knowledge. However, mutual knowledge implies nothing about what people know about other people's knowledge. 

Common knowledge is a related but stronger notion. A fact is called \textit{common knowledge} if all the players know the fact, and they also know that everyone knows the fact, and they also know that everyone knows that everyone knows the fact, and so on ad infinitum. So by making a public announcement at the beginning, Tanya made the fact that there was at least one red hat common knowledge instead of just mutual knowledge.

By definition, we know that any fact that is common knowledge is also mutual knowledge, but not the other way around.

The puzzle above is known under many disguises: cheating wives, cheating husbands, professors and their published mistakes, blue-eyed people on an island, and muddy children \cite{F, MG} are most popular among them. The concept of common knowledge was introduced in philosophical literature in 1969 \cite{F,FHMV}. 

In this paper we play with common knowledge and introduce a lot of new common knowledge puzzles. We start by defining circular and simultaneous games in Section~\ref{sec:prel}. In Section~\ref{sec:sight}, we introduce variations of this puzzle where the sages have various sight problems. What will happen to the game if one of the sages is blind? We also have puzzles where sages are near-sighted or far-sighted.

In Section~\ref{sec:numbers} we replace colored hats with integers. The colors can be matched to numbers, or the numbers can be written on the hats instead of colors, or the numbers can be attached on Post-It notes to sages' foreheads. Introducing numbers is more than just introducing many hat colors. We will actually use the values of the integers. 

In the paper \textit{A Headache Causing Problem} by Conway, Paterson and Moscow \cite{CPM} the notes with numbers were nailed to the foreheads of the players thus causing a headache. Our evil emperor is cruel, but for our purposes it doesn't matter how numbers are attached. The important fact is that the sages see other people's numbers, but not their own.

In Section~\ref{sec:md}, we solve puzzles where the emperor announces the maximum difference between the numbers on the sages' foreheads. In Section~\ref{sec:SorP}, the emperor announces a number which is either the sum or the product of all the numbers on the foreheads.

We discuss solutions to all the puzzles. But for some of the puzzles we invite the readers to solve them first. In those cases the solutions are delayed to the end of the paper in Section~\ref{sec:solutions}.

\section{Circular and Simultaneous Games. Preliminaries.}\label{sec:prel}

We begin by establishing some rules for the puzzles we will discuss. We assume that there are $N$ players, and each of them can see everyone else, unless otherwise noted. In most puzzles the players sit in a circle. We also assume that the players are infinitely intelligent and understand common knowledge. 

We also assume that each player says either YES or NO. They say YES if they know what is on their forehead, and NO if they do not. For example, saying YES is equivalent to saying ``I know my hat color'' or ``I know my number.'' 

Let's go back to Puzzle~\ref{puzzle:intro}. We now assume that all the three sages say YES or NO at the same time. 
Then they wait the agreed amount of time and proceed to the next announcement. We call each time step a \text{round} of the game. For example, if there are two red hats, then everyone says NO in the first round. In the second round two people with red hats say YES, and the person with a blue hat says NO. In the third round, everyone says YES. By the way, if a person says YES in a round, then s/he knows the hat color and says YES in all the next rounds if any.

We call the game where all the sages speak together \textit{a simultaneous game}. 

Consider this generalization of Puzzle~\ref{puzzle:intro}.

\begin{puzzle}[A generalization of the intro puzzle.]
An emperor puts his $N$ sages to a test. He shows them an infinite supply of blue and red hats. The sages close their eyes and the emperor puts a hat on each of their heads. They open their eyes and see the hats on other sages but not on themselves. The emperor says that at least one of them has a red hat. If at least one of them can figure out the color of their hat, they all can go free. The sages are playing a simultaneous game. What will happen?
\end{puzzle}

In this puzzle, all of the people with the same hat color see the same set of colors, so they all behave in the same manner: they will say YES or NO at the same time. Similar to the introduction game, if on round $k$ everyone says NO, it becomes common knowledge that there are at least $k+1$ red hats. Thus every sage with a red hat will say YES for the first time during round $m$, where $m$ is the number of red hats. Then in round $m+1$, all sages with blue hats will know their color and say YES.

We also want to introduce another variation of this game, where we assume that the players announce their knowledge in a circular manner one after the other, starting with the first player. We call it \textit{a circular game}. In a circular game a \textit{turn} is for one person to announce their knowledge. The next person in a circle can use this information to update their knowledge. When all $N$ people announce their knowledge and finish the circle, we call it a \textit{round}. In this case, the round consists of $N$ turns.

Now we analyze the same puzzle if it were played as a circular game.

\begin{puzzle}[A circular version of the intro puzzle.]
An emperor puts his $N$ sages to a test. He shows them an infinite supply of blue and red hats. The sages close their eyes and the emperor puts a hat on each of their heads. They open their eyes and see the hats on other sages but not on themselves. The emperor says that at least one of them has a red hat. If at least one of them can figure out the color of their hat, they all can go free. The sages are playing a circular game. What will happen?
\end{puzzle}

This game is known and was covered in \cite{FHMV}. We want to explain what is going on as we have similar problems in other circular games later. 

\textbf{Answer.} In a circular game, during the first round the last red-hatted sage to speak figures out that he has a red hat, and everyone after him in a circle knows that they have a blue hat. No one else will ever figure anything out, even after the first round.

\textbf{Explanation.} 
Consider the first person, Kevin. If Kevin sees only blue hats, he knows his hat is red. But if he 
sees a red hat, then he cannot determine his own hat's color. Thus, if the first person says he knows his hat color, then everyone else will know immediately that their hats are all blue. 

Suppose Kevin sees a red hat and announces that he does not know his hat color. Everyone else now has common knowledge that there must be a red hat among the sages excluding Kevin. If the second person, Armaan, sees only blue hats among everyone else besides the first person, he knows his own hat must be red. 

In this manner, as long as a student announces that he doesn't know his hat's colors, everyone knows that there must be a red hat between that person and the end of the circle. As soon as someone only sees blue hats after them in the circle, they know their hat is red. As soon as someone announces that they know their color, everyone else who has yet to speak knows their hat is blue.

In general, the order in which players speak affects who will be able to figure out their hat colors. For example, if there is only one red hat and the person with the red hat is first, then everyone will figure out their color. If the red-hatted person is the last to speak, then only they will know their hat color. 

In this puzzle, the number of people who can figure out their color varies between 1 and the number of blue hats plus 1 depending on the order in which they speak. In general, the circular game is more complicated than the corresponding simultaneous game.

\section{The Mind can Know what the Eye doesn't See}\label{sec:sight}

\subsection{A Blind Person.}

We now want to add a blind person, Alice, to our group of infinitely intelligent sages sitting in a circle. 

\begin{puzzle}[Blind Alice and the quest for the red hat.]
The total number of people is $N$ and they each have a red or a blue hat on their heads, which is common knowledge. Additionally, it is common knowledge that Alice is the only blind person and that at least one person has a red hat. What will happen?
\end{puzzle}

We assume that Alice knows how many people there are. For a circular game, we also assume that Alice knows where she sits and when it is her turn to speak.

The question is what can they figure out about their hats.

\textbf{Simultaneous game.} 

Suppose Alice has a blue hat. This case is nearly identical to the original hats problem, where a time step is used. Let $r$ be the number of red hats. During round $r$, all the red-hatted people will say YES, and since Alice has a blue hat, her inability to see does not matter up to this point. After that Alice and all of the other people with blue hats know that they have blue hats, so they say YES during round $r+1$.

Suppose Alice has a red hat. Here, the other people cannot make conclusions like they could in the original hats problem, since the blind person cannot see how many other red hats there are, resulting in no one else being able to distinguish between whether they have red or blue. The sages have lost their starting point where only 1 person having red knows the color. If Alice has the only red, she will not say YES in the first round. However, Alice knows that if she has a blue hat then by round $N-1$ someone would say YES. This means that after $N-1$ rounds, Alice will know that she must have a red hat, so she will say YES on round $N$. No one else will be able to say anything more. 

\textbf{Circular game.} 

We can apply the logic from the circular hat game without blind people up until Alice speaks. 

1) Suppose that Max, who has a turn before Alice, says YES in the first round. That means Max is the last person with a red hat and everyone else has a blue hat. It doesn't matter that Alice is blind. The game proceeds the same way as if she is not blind.

2) Suppose Alice is the last person and no one before her figures out their hat color. By the same logic as in the game without a blind person, Alice knows that she must have a red hat. No other sage can figure out their colors.

3) Now suppose that Alice is not the last person in the circle, and no one before her said YES. That necessitates that either Alice or someone after her has a red hat, so Alice is bound to say NO.

3a) If Alice has a blue hat, people in the circle after Alice know that there is at least one red hat among them. This is similar to what we discussed before. The last person in the circle with a red hat, call him Kunal, will know that he has a red hat and will say YES. Everyone after Kunal will realize they have blue hats. On the next round Alice knows that she has a blue hat, but no one else can figure out their hat colors.

3b) If Alice has a red hat, other people after her do not know if any of them have a red hat, so everyone says NO. On the next round, Alice knows her hat color, and she is the only person who can eventually figure out the hat color.

Ironically, in both versions of the game, the blind sage is the only sage who can always figure out their hat color.

\subsection{Near-sighted People}

We now create a variation of the puzzle where all of the players are \textit{near-sighted}. This means that they can only see their closest neighbors in the circle.

\begin{puzzle}[Near-sighted people.]\label{puzzle:ns-sim}
There are only red and blue hats. There are $N$ near-sighted sages with hats sitting in a circle playing a simultaneous game. The emperor announces to them that there is exactly one red hat. Not all of the sages eventually figured out their colors. What is $N$?
\end{puzzle}

We leave it to the reader to think about this puzzle. The solution is in Section~\ref{sec:solutions}. Meanwhile we will analyze the circular game with $N$ near-sighted people.

\textbf{Circular game.}

Here is the summary of what happens, where we number the people from 1 to $N$ and by $r$ we denote the index of the person with the red hat:

\begin{itemize}
\item If $N=3$, everyone sees all the other hats and can deduce their own hat color immediately.
\item If $N > 3$ and $r = 2$ or $r = N$, the first person will say YES because they see a red hat, so everybody else will know either the second or last person has the red hat. The second person will say NO, then everybody else after him up to the last persons will say YES. The last person will say NO. This exact same sequence of answers will occur in later rounds as well, so, in this case, everybody except the second and last people will eventually figure out their colors.
\item If $r$ is odd and $N>3$ is even, the first person doesn't see a red hat and says NO. It becomes common knowledge that the second and the last person have blue hats. The second person will say YES. This statement doesn't provide any new information so the next sage says NO and this continues. Odd-numbered people say NO because they see both blue hats, which allows the next even-numbered person to say YES without providing any new information to the person after that. Only the even-numbered people will eventually figure out their colors.
\item If $r$ is even and $2 < r < N$, the first person will say NO, the second person will say YES,
and it will continue to alternate until sage number $r-1$. Sage $r-1$ will say YES, breaking the alternating pattern.
Then everybody knows that sage number $r$ has the red hat, which allows everyone to figure out their hat color.
\item If $r$ is odd and $N$ is odd, the first person doesn't see a red hat and says NO. It becomes common knowledge that the second and the last person have blue hats. The second person will say YES. This statement doesn't provide any knew information so the next sage says NO and this continues. Odd-numbered people say NO because they see both blue hats, which allows the next even-numbered person to say YES without providing any new information to the person after that. This conitnues until the very end, when the last person knows their hat color because the first sage said NO. Eventually the even-numbered people and the last person will say YES.
\end{itemize}

\subsection{Far-sighted People}

For symmetry, we want to study a game where all the players are \textit{far-sighted people}, which means that they can see everyone but their neighbors in the circle.

\begin{puzzle}[Ten far-sighted people.]
There are only red and blue hats. There are $10$ far-sighted people playing a simultaneous game. The emperor announces to them that there are exactly three red hats. How should the emperor place the hats to minimize the number of sages that will figure out their hat colors?
\end{puzzle}

In the first round, if someone sees all three red hats, then they know their hat is blue. Also, if someone sees all seven blue hats, they know their hat is red.

There are seven people with blue hats, so by the Pigeonhole Principle, at least three of the people with blue hats must be in a row. Therefore, the person in the middle of these three people is guaranteed to see three red hats, so in the first round, at least one person will say YES.

Now we consider cases depending on how the blue hats are split up by the red hats. The different ways that the blue hats can be split up are as follows: $(7, 0, 0)$, $(6, 1, 0)$, $(5, 1, 1)$, $(5, 2, 0)$, $(4, 2, 1)$, $(4, 3, 0)$, $(3, 2, 2)$, and $(3, 3, 1)$. Because the players are sitting in a circle, we can assume that the circle is rotated in such a way as to make the largest group of consecutive blue hats at the end. When we write the possible orderings, we denote a red hat by R and a blue hat by B. We will leave the last group of $B$s out. For example, the $(7, 0, 0)$-split corresponds to the hat ordering RRR, which means sages 1 through 3 have red hats.

In Table~\ref{table:10fsp} we show the sages who say YES in the first round in each of the cases.

\begin{table}[h!]
\centering
\begin{tabular}{|c | c | c | c|} 
 \hline
 Case \# & Split & Hats & YESes \\ 
 \hline\hline
 1 & $(7,0,0)$ & RRR & 2, 5, 6, 7, 8, 9 \\ 
 2 & $(6,1,0)$ & RBRR or RRBR & 6, 7, 8, 9 \\
 3 & $(5,2,0)$ & RRBBR or RBBRR & 7, 8, 9 \\
 4 & $(5,1,1)$ & RBRBR & 7, 8, 9 \\
 5 & $(4,3,0)$ & RRBBBR or RBBBRR & 4, 8, 9; or 3, 8, 9 \\
 6 & $(4,2,1)$ & RBRBBR or RBBRBR & 8, 9 \\
 7 & $(3,2,2)$ & RBBRBBR & 9 \\
 8 & $(3,3,1)$ & RBRBBBR or RBBBRBR & 5, 9; or 3, 9 \\
 \hline
\end{tabular}
\caption{The first round for the game with 10 far-sighted people}
\label{table:10fsp}
\end{table}

We see that in cases 1, 5, 7, and 8, the configuration is uniquely determined by the list of sages who say YES in the first round. This means that in the second round, everyone else will say YES.

Consider case 2 with a $(6,1,0)$-split corresponding to RBRR or RRBR hat orderings. These are the only hat orderings where 4 consecutive sages say YES in the first round. After the first round, it becomes common knowledge that people number 5 through 10 are blue and people number 1 and 4 are red. These people will say YES in the second round. But people number 2 and 3 have no way of knowing who between them has the red hat. In this case the game stabilizes after two rounds, and only 8 people will ever know their color.

Consider cases 3 and 4 with a $(5,2,0)$-split corresponding to RRBBR or RBBRR hat orderings, or a $(5,1,1)$-split corresponding to RBRBR hat ordering. These are the only hat orderings where 3 consecutive sages say YES in the first round.  After the first round, it becomes common knowledge that people number 6 through 10 have blue hats and people number 1 and 5 have red hats. In the case RRBBR, the sage number 4 sees a red hat on person 2, and can deduce that they have a blue hat. Hence, in the second round, people number 1, 4, 5, 6, and 10 say YES. Similarly, in the case  RBBRR, the second person realizes that they have a blue hat. Hence, in the second round sages number 1, 2, 5, 6, and 10 say YES. The second round disambiguates the configurations. Therefore, in the third round everyone knows their hat color.

Consider case 6 with a $(4,2,1)$-split corresponding to RBRBBR or RBBRBR hat orderings. These are the only hat orderings where 2 consecutive sages say YES in the first round.  After the first round, it becomes common knowledge that people number 2, 5, 7, 8, 9, and 10 have blue hats and people number 1 and 6 have red hats. They all will say YES in the second round. After that, people number 3 and 4 can't figure out who between them is red. The game stabilizes after the second round and only 8 people can figure out their colors. 

Therefore, at least eight people will always be able to figure out their hat colors, and to achieve this minimum, the emperor should place the hats in one of the following configurations: $RBRR$, $RRBR$, $RBRBBR$, or $RBBRBR$.

Here is another puzzle on the same topic for our readers to solve. The solution is in Section~\ref{sec:solutions}.

\begin{puzzle}[Nine far-sighted people.]\label{puzzle:9fs}
There are only red and blue hats. There are $9$ far-sighted people playing a simultaneous game. The emperor announces to them that there are exactly three red hats. In the first round everyone says NO. What will they say in the second round?
\end{puzzle}

\section{From Colors to Numbers}\label{sec:numbers}

This time instead of hats, we have integers placed on each players' foreheads. The integers can be written on Post-It notes, or we can reuse hats and write numbers on hats, one integer per person. As before, participants see other  people's numbers but not their own. Here people can say ``I know my number,'' or ``I don't know my number.'' We will often shorten this to YES for the former and NO for the latter.

\begin{puzzle}[Zeros and Ones.]
Every player has either 0 or 1 are on their foreheads, and this is common knowledge. The players need to guess their numbers in a simultaneous or a circular game. The emperor announces that there is at least one zero. How will the game go?
\end{puzzle}

This game is equivalent to the game with hats. We can replace each zero with a red hat and each one with a blue hat.

We can use the game with hats to understand many other games with numbers. For example, we do not need to restrain ourselves to only two sages and only zeros and ones. Suppose we have $N$ sages with numbers on their foreheads. Suppose the emperor announces that there exists at least one person with a given number $x$. Then we can replace all of the notes with $x$'s with red hats and all of the other numbers with blue hats. Then, by our previous analysis, the sages will figure out their hat colors. That means, the people with $x$'s will figure out their number, but the other people will not.

We leave the next puzzle as an exercise to the reader.

\begin{puzzle}[Ten near-sighted people in a line.]\label{puzzle:10ns}
Ten near-sighted sages are sitting in a line. The emperor announces that they each have a positive integer on their foreheads and the sum of their numbers is either 30, 31, or 32. The leftmost person, Aidan, starts by saying whether or not he knows his own number. Then Josh, Aidan's neighbor to the right says whether he knows his number. If not all the people say YES in the first pass, they continue starting from Aidan again. The sages are lucky: the emperor assigned numbers in such a way that they all can figure out their integers in the smallest number of turns. What are the numbers?
\end{puzzle}

The solution is in Section~\ref{sec:solutions}.

\section{Maximum Difference}\label{sec:md}

The following puzzles and variations of it appear a lot on the Internet and in different books.

\begin{puzzle}[Two consecutive numbers.]\label{puzzle:2cons}
Alice and Bob are told that they have two consecutive non-negative integers on their foreheads. Alice says, ``I don't know my number.'' Then Bob says, ``I also don't know my number.'' Then Alice says, ``I still don't know my number.'' And Bob counters, ``Now I know my number.'' What are their numbers if Bob's number is odd?
\end{puzzle}

The solution is in Section~\ref{sec:solutions}.

In this section, we discuss a variety of puzzles where our evil emperor announces the maximum difference of the numbers on the player's foreheads. 

\textbf{General Setting.} There are $N$ people who are told that they have non-negative integers on their foreheads. The emperor announces that the largest difference between their numbers is $D$. We denote the maximum number on the sages' foreheads as $M$ and the minimum as $m$. Let the number of people with the maximum be noted as $P$. We call the people who do not have $m$ or $M$ \textit{normal}. We call the people who do have $m$ or $M$ \textit{extreme}.

The puzzle at the beginning of this section is a particular example of the maximum difference puzzle corresponding to a circular game in which $N=2$ and $D=1$.

\subsection{Two people}

Let our players be named Alice and Bob, as the most famous characters in cryptography. It is very convenient that their names start with $A$ and $B$. So, we denote the numbers on their foreheads as $A$ and $B$ correspondingly. Consider a circular game, where Alice starts, as she is earlier in the alphabet. Alice can only know her number if $B < D$. If on the first turn Alice says NO, then Bob knows that $B \geq D$. That means, he can figure out his number only if $A-D < D$, equivalently $A < 2D$, which implies that he must have $A+D$. That means, he can only figure out his number if $B < 3D$. If Bob says NO, then Alice doesn't have a number between 0 and $2D-1$. Continuing in a similar manner, and we see that during round $t$, after $t-1$ NOs, the person knows that their number must be more than $(t-1)D -1$.

Suppose the maximum is $M$. The person with the maximum sees $M-D$ and is choosing between $M$ and $M-2D$. So this person guesses their number as soon as $(t-1)D-1\geq M-2D$. This will occur on turn $\left \lceil \frac{M+1}{D} -1\right \rceil$. In particular, if Alice has the maximum, she will say YES in round $k$ and turn $2k-1$, where $k = \left \lceil \frac{M+1}{2D}\right \rceil$. If Bob has the maximum, he will say YES in round $k$ and turn $2k$, where $k = \left \lceil \frac{M+1}{2D}-\frac{1}{2}\right \rceil$.

We now consider a simultaneous game. In the first round, a person knows their number if and only if the number they see is less than $D$. If both say that they don't know their number, then it becomes common knowledge that their numbers are not less than $D$. Thus if they see a number in the range $[D,2D-1]$, they know their number in the second round. This continues, and if on the $k$-th round they both say they don't know their number, then they know that their numbers are more than $kD-1$. Continuing, we see that the person with the largest number $M$, will say YES on round $k$, where $M \in [(k-1)D,kD-1]$.

In both games, the person with the larger number will know their number first, and the person with the smaller number will know their number right after that.

\subsection{$N$ people. Consecutive numbers}

In this section we study a game where the integers are non-negative, consecutive and distinct, and this fact is common knowledge. It follows that $N=D+1$. Normal people immediately know their numbers. The person with the maximum, say Kunal, knows his number before the start of the game if and only if he sees a zero. We assume that $N > 2$, as we discussed the case of two people above. This means that normal people exist.

\textbf{Game 1: simultaneous.}  We have two cases. 

Case 1. $M=D$. This means $m=0$. Then during the first round all normal people and the person with the maximum say YES. Then, in the second round, the person with the minimum will realize, that they have 0, so they will say YES.

Case 2. $M>D$. During the first round, all normal people say YES. In the second round, the two extreme people can deduce what their numbers are.

\textbf{Game 2: circular.} Note that, as in the simultaneous case, the extreme people always know that they are extreme and the normal people always say YES. Therefore, if an extreme person says NO, everyone realizes that this person is extreme and so the second extreme person will know that they are extreme and say YES.

Suppose $N =3$ and there are three people total: Alice, Bob, and Cathy. Alice goes first followed by Bob.

If Alice says NO, then everyone else knows that Alice is extreme and can figure out their numbers. Alice doesn't have a way to figure out her number.

If Alice says YES and she doesn't have a 2, she must be normal, and both Bob and Cathy will know their numbers after that. 

If Alice says YES and she has a 2, then the the set of numbers could only be $(0, 1, 2)$ or $(1, 2, 3)$. In any case, the person with 1 will know their number after Alice's turn, and the other person will never know their number.

Here is another puzzle left for the reader. 

\begin{puzzle}[Four people with consecutive numbers.]
Four people, Aaron, Bibi, Cal, and Dorothy are playing a circular game of maximum difference. It is common knowledge that their numbers are non-negative and consecutive. Everyone says YES in the first round. What can we conclude about the numbers and the order in which they are placed?
\end{puzzle}

\textbf{Solution.} We know that Aaron says YES if he is normal or if he has the maximum and it is 3. 

\textbf{Case 1. Aaron doesn't have a 3.} Then it becomes common knowledge that he is normal. If his number is $a$, then it is common knowledge that numbers $a-1$, $a$, and $a+1$ are present on player's foreheads. This means it doesn't matter in which order people with numbers $a-1$, $a$, and $a+1$ go: they do not provide any new information for the fourth person who has either $a-2$ or $a+2$. The fourth person is extreme. The only way for this person to know their number is if they knew it from the start. This could only happen if this number is 3 and $a=1$. 

\textbf{Case 2. Aaron has a 3.} That means the set of numbers must be one of the following: $\{0,1,2,3\}$, $\{1,2,3,4\}$, and $\{2,3,4,5\}$. It becomes common knowledge that 2 is one of the numbers. That means, when the person with 2 speaks, there is no new information provided. Let's ignore the person with 2 and consider the next person. This person can't be extreme, as the extreme person can't figure out their number yet: both extreme 0 and 4 see the same $\{1,2,3\}$ and both extreme 1 and 5 see the same $\{2,3,4\}$. That means the next person excluding 2 must be normal. If it is 1, then the last person knows that there is a zero. If it is 4, then the last person knows that there is a 5.

Overall, everyone says YES in three cases: 
\begin{itemize}
\item The numbers are $\{0,1,2,3\}$ and Aaron has 1.
\item The numbers are $\{0,1,2,3\}$, Aaron has 3 and the person with 1 goes before the person with 0.
\item The numbers are $\{2,3,4,5\}$, Aaron has 3 and the person with 4 goes before the person with 5.
\end{itemize}

\subsection{$N$ people. Not necessarily distinct, $D=1$, simultaneous.} 

In the first round, a person knows their number if and only if everyone else has 0. This means that $P = M = 1$. In this case, the person with 1 knows their number in the first round, and everyone else will know their number in the second round. 

Suppose that no one knows their number in the first round. This means that it is common knowledge that we cannot have $P=M=1$. If someone knows their number on the second round, then they must have value 1; moreover, there must be exactly two such players, and the rest must have value 0. If two people say YES in the second round, then the rest will know their numbers in the third round. If no one says YES in the second round it becomes common knowledge that we cannot have $(M, P)$ in $\{(1,1), (1,2)\}$. In the next round people who see two 1s and a bunch of zeros can conclude that they have a 1.

Continuing in this manner, we see that if there are $r$ ones and $N-r$ zeros, all $r$ people with ones will say YES in round $r$. In the next round, the people with zeros will know they have zero and say YES. If no one said YES in round $r$, it becomes common knowledge that we cannot have $M=1$ and $P \leq r$. On round $N$ everyone will realize that there are no zeros. If someone sees $N-1$ ones, they will know their number has to be 2.

If someone sees $N-2$ ones in round $N+1$, they know their number has to be 2, and so on until round $2(N-1)$, when everyone will figure out that there are no ones. During this round if someone sees $N-1$ twos, they will know their number must be 3. This process will continue until people with the maximum will figure out their numbers in round $m(N-1)+P$. People with the minimum will figure out their number in the next round.

\subsection{$N$ people. Not necessarily distinct, $D>1$, simultaneous.} 

For a final touch in this section on maximum difference, we leave the readers with a crowning puzzle to solve.

\begin{puzzle}[One, Zero, Two, Three.]\label{puzzle:mdsummary}
Several sages are playing the maximum difference simultaneous game. In the first round, one person knows their number. In the second round, nothing new happens. In the third round, two more people know their numbers. In the fourth round, three more people know their numbers. How many people are there are what are the numbers?
\end{puzzle}

The solution is in Section~\ref{sec:solutions}.

\section{Sum or product}\label{sec:SorP}

This section is devoted to puzzles where the evil emperor announces that an integer is either the sum or the product of all the numbers on the foreheads of the players.

\textbf{General setting.}
Suppose that a finite number $N$ of sages sit in a circle.
The emperor puts Post-Its on everyone's forehead, with integers written on each
of them, and announces that a certain number $M$ is either the product or the
sum of all of the numbers. It is common knowledge among all of the sages that all of their numbers are positive integers. People see each other's numbers, but do not see their own numbers.

First, we discuss scenarios where one sage, say Alice, immediately knows her number. We have three cases.

\begin{itemize}
\item \textbf{It-must-be-the-product.} The numbers that Alice sees sum up to at least $M$.
\item \textbf{It-must-be-the-sum.} The numbers that Alice sees have a product that doesn't divide $M$.
\item \textbf{Whatever.} Whether $M$ is the sum of the product Alice gets the same result for her number. 
\end{itemize}

Keep in mind: Alice has at most two possibilities for the integer on her forehead and as soon as she knows whether $M$ is the sum or the product she knows her number.

\subsection{Sum or product: two people.}

The following known folklore puzzle inspired the research in this section.

\begin{puzzle}[50/50 Sum or product. Circular.]\label{puzzle:50circ}
Alice and Bob have positive integers on their foreheads. They see each other's numbers, but do not know their own number. It is common knowledge that they have positive integers. An observer announces to Alice and Bob that number 50 is either the sum or the product of their numbers. After that, Alice says to Bob, ``I do not know my number,'' to which Bob replies, ``I do not know my number either.''

What is Bob's number?
\end{puzzle}

Here is another puzzle which is a simultaneous version of the puzzle above.

\begin{puzzle}[50/50 Sum or product. Consecutive.]\label{puzzle:50sim}
Alice and Bob again have positive integers on their foreheads. Again, they see each other's numbers, but do not know their own number, and it is common knowledge that they have positive integers. An observer announces to Alice and Bob that number 50 is either the product or the sum of their numbers. After that, Alice and Bob both say to each other, at the same time, ``I don't know my number'' twice. What are their numbers?
\end{puzzle}

The solutions to both puzzles are in Section~\ref{sec:solutions}.

We provide a complete analysis of the sum-or-product simultaneous and circular games for any announced integer $M$ and two people. We denote Alice's number by $A$ and Bob's number by $B$.

The three cases when Alice knows her number from the start for this particular set up are as follows.

\begin{itemize}
\item \textbf{It-must-be-the-product.} This could only happen if $M=B$.
\item \textbf{It-must-be-the-sum.} This could only happen if $B$ doesn't divide $M$.
\item \textbf{Whatever.} This could only happen if $M =4$ and $B=2$.
\end{itemize}

Thus a sage doesn't know their number immediately if they see a proper divisor and it is not the case of $M=4$ and $A=B=2$.

Back to the game. We start with the cases when they both know their numbers at the very beginning and, hence, both say YES in the first round of either a circular or a simultaneous game.

\begin{itemize}
\item \textbf{It-must-be-the-product.} This is the case when $M =A=B=1$.
\item \textbf{It-must-be-the-sum.} Both numbers are non-divisors of $M$.
\item \textbf{Whatever.} This could only happen if $M =4$ and $A=B=2$.
\end{itemize}

Now we exclude the cases above and see whether Alice's first statement is useful to Bob. Bob has two choices for his number $M-A$ or $\frac{M}{A}$ and Alice's statement helps if it can differentiate between these two possibilities.

Alice would say YES for both $M-A$ and $\frac{M}{A}$ only if $A$ is 1 and $M > 2$. Alice would say NO for both $M-A$ and $\frac{M}{A}$ if $M-A$ is a divisor of $M$. This could only happen when $A = \frac{M}{2}$. In all other cases, that is, when $M-A$ is a non-divisor and $A >1$ or if $M=2$ and $A=1$, Alice's statement will eliminate one of the possibilities for Bob's number.

Now we consider cases where one of the numbers is 1. Recall that we already excluded the case $M=1$.

\begin{itemize}
\item $A=B=1$ and $M=2$. In the simultaneous game both say NO in the first round, then both say YES in the second. In the circular game, Alice says NO, then Bob says YES, and Alice will never know whether she has 1 or 2.
\item $M=2$ and it is the product. In the simultaneous game, the person who sees 2 says YES in the first round while the other person says NO. The other person will realize that their number is 2 and say YES in the second round. In the circular game, if Alice sees 2, she says YES, and Bob realizes that he must have 2 and also says YES. If Alice sees 1, she says NO, then Bob sees 2 and says YES, and Alice will never know whether she has 1 or 2.
\item $A$ or $B$ equals 1 and $M > 2$. In this case, in both simultaneous and circular games, the person with 1 will immediately know their number and the other person will not know their number ever. 
\end{itemize}

Now supposes $A > 1$, and $B > 1$. Excluding the cases we discussed, we can assume that $M > 4$.

\begin{itemize}
\item One of $A$ and $B$ is a non-divisor of $M$ and the other is a divisor that is larger than 1. In the simultaneous game, the person who sees a non-divisor says YES in the first round and the other person says NO. In the second round, the person who sees a divisor realizes that $M$ must be the sum and says YES. In the circular game, if Alice sees a non-divisor she says YES, and Bob realizes that $M$ must be the sum and also says YES. If Alice sees a divisor, she says NO. Bob says YES, and Alice will never know her number.
\item Both $A$ and $B$ are divisors where neither is equal to 1 or $\frac{M}{2}$. In the simultaneous game, both say NO in the first round, then both say YES in the second round. In the circular game, as Alice says NO in the first turn, Bob knows he must have a proper divisor. In this case, $M - B$ is not a proper divisor, so Bob knows his number and says YES. After that Alice will never know her number.
\item $A=B = \frac{M}{2}$. In the simultaneous game, before the first round, both of them don't know whether they have 2 or $\frac{M}{2}$. Both say NO in the first round. In the second round, both people know that if they had $2$, then the other person who sees a 2, thinks they have either $\frac{M}2$ or $M-2$. But $M-2$ triggers YES in the first round. Thus if someone sees a 2, they would say YES in round 2. In actuality, they both say NO in round 2, which makes them realize that neither of them see a 2. Thus in the third round, they both know that they  have $\frac{M}{2}$ and both say YES. In the circular game, Alice says NO in the first round. In this case, Bob knows that he has either $\frac{M}{2}$ or 2 and Alice's first statement doesn't help. He says NO on his turn.  On the third turn, Alice realizes that she has $\frac{M}{2}$, as otherwise Bob would have guessed his number. But Bob will never know whether he has 2 or $\frac{M}{2}$.
\item Alice and Bob have $\frac{M}{2}$ and 2 as their numbers. In the simultaneous game, suppose Bob has 2. Both say NO in the first round. For the second round, Alice can deduce that she can't have $M-2$ due to the fact that if she did, Bob would have said YES in the first round as $M-2$ doesn't divide $M$. Alice realizes that she must have $\frac{M}{2}$ and says YES. Bob can't differentiate between 2 and $\frac{M}{2}$ and says NO. In the third, round Bob knows if he had $\frac{M}{2}$, Alice would not have said YES by the previous case. Thus Bob knows he has 2 and says YES in the third round. In the circular game, suppose Bob has 2. Then Alice starts by saying NO. Bob also says NO. On the third turn, Alice realizes that she can't have $M-2$, and says YES. But Bob will never know whether he has 2 or $\frac{M}{2}$. If Alice has 2, then again Alice starts by saying NO. Bob realizes that he can't have $M-2$ and says YES. Alice realizes that she can't have $\frac{M}{2}$ as Bob would have said NO. She says YES on her turn.
\end{itemize}

We summarize the answers here to show how the game goes. First we describe the simultaneous game. The first answer in each round is Alice's. We use a semicolon to separate the rounds. We show the last new round that is different from the next one.

\begin{itemize}
\item YES, YES: $A=B=M=1$; or neither $A$ nor $B$ divide $M$; or $M=4$ and $A=B=2$.
\item YES, NO: $A=1$, $M >2$.
\item NO, YES: $B=1$, $M >2$.
\item NO, NO; YES, YES: $A=B=1$ and $M=2$; Both $A$ and $B$ are divisors where neither is equal to 1 or $\frac{M}{2}$.
\item YES, NO; YES, YES: $A=1$, and $B=M=2$; $A>1$ is a divisor and $B$ is a non-divisor;
\item NO, YES; YES, YES: $B=1$ and $A=M=2$; $B>1$ is a divisor and $A$ is a non-divisor.
\item NO, NO; NO, NO; YES, YES: $M > 2$, $M \neq 4$ and $A=B = \frac{M}{2}$.
\item NO, NO; YES, NO; YES, YES: $A=\frac{M}{2}$ and $B=2$.
\item NO, NO; NO, YES; YES, YES: $A=2$ and $B=\frac{M}{2}$.
\end{itemize}

Now we summarize the circular game.

\begin{itemize}
\item YES, YES: $A=B=M=1$; or $M=4$ and $A=B=2$; or or $A=1$ and $B=M=2$; $B$ is a non-divisor.
\item YES, NO: $A=1$, $M >2$.
\item NO, YES: $B=1$, $M >2$; $B=1$ and $M=2$; $B>1$ is a divisor and $A$ is a non-divisor; Both $A$ and $B$ are divisors where neither is equal to 1 or $\frac{M}{2}$.
\item NO, NO; YES, NO: $M > 4$, and $A= \frac{M}{2}$.
\item NO, YES; YES, YES: $M > 4$, $A=2$ and $B=\frac{M}{2}$.
\end{itemize}

We notice that in any case at least one of them always figures out their number. In the simultaneous game, only one person knows their number at the end only when $M > 2$ and one of $A$ or $B$ is 1. In a circular game only one person ends up knowing their number when:

\begin{itemize}
\item $M > 2$ and one of $A$ or $B$ is 1.
\item $B=1$ and $M=2$; $B>1$ is a divisor and $A$ is a non-divisor; Both $A$ and $B$ are divisors where neither is equal to 1 or $\frac{M}{2}$.
\item $M > 4$, and $A= \frac{M}{2}$.
\end{itemize}

\subsection{$N$ people in a simultaneous game. The announced number is prime}

We now have $N$ people and the announced number $M$ is prime. Here are the three cases when Alice knows her number from the start.

\begin{itemize}
\item \textbf{It-must-be-the-product.} This could only happen if Alice sees all ones and $N > M$, or if she sees all ones and a single $M$.
\item \textbf{It-must-be-the-sum.} This could only happen if Alice sees at least one number that is not equal to 1 or $M$.
\item \textbf{Whatever.} This could never happen.
\end{itemize}

Alternatively, the only situation when Alice doesn't know her number is when she sees all ones and $N < M+1$.

First we consider a simultaneous game.

\begin{itemize}
    \item \textbf{There are at least two non-ones.} If there are at least two non-ones, then every person will see at least one non-one. Since there are two non-ones in total, none of them is $M$. Since every person sees one non-divisor of $M$, they know $M$ is the sum, and then they can figure out their number. Therefore, in the first round, everyone will know their number.
    \item \textbf{There is exactly one non-one.} If there is exactly one non-one, either the non-one is $M$ and $M$ is the product, or the non-one is $M-N+1$ and $M$ is the sum. If the non-one is $M$, all of the other people will see $M$, so they know it must be the product and they must have 1. If the non-one is $M-N+1$, everyone else will see a non-divisor of $M$, so they will know that $M$ is the sum and will be able to figure out their numbers. However, in both cases, the person with the non-one will never know what they have.
    \item \textbf{All the numbers are one.} In this case, $M$ is the sum and $M = N$. Everybody sees $N-1$ ones, so they do not know if they have $1$ or $M$. Let Alice be one of the people. Then when Alice hears that no one else knows their number, she knows that she has one because if she had $M$, everyone else would have known that they had one. The same holds for all of the other people. Therefore, no one will know their number in the first round, but everyone will know their number in the second round. If $M=1$, everyone knows their number immediately.
\end{itemize}

Now we proceed to a circular game with the same cases. We assume that Alice is the first person.

\begin{itemize}
    \item \textbf{There are at least two non-ones.} As we explained in the simultaneous case everyone knows their numbers before the game begins. Thus, everyone says YES.
    \item \textbf{There is exactly one non-one.} Every one who has a one know their number before the game starts. So they all say YES. The person with a non-one has no way of figuring out whether they have $M$ or $M-N+1$ and $M$. So this person will never know their number.
    \item \textbf{All the numbers are ones.} If $M=N > 1$, then Alice doesn't know if she has $M$ or 1. She says NO. Everyone else concludes that they must have ones and say YES. Alice will never know her number. If $M=1$ everyone knows their numbers.
\end{itemize}

\subsection{$N>2$ people in a simultaneous game. The announced number is semiprime}

Currently, we have results when $N=2$ people, and also when the announced number $M$ is prime. What if the number is semiprime and we have more than two people? Here are some puzzles.

\begin{puzzle}[Everyone Says NO.]\label{puzzle:35}
Three people are playing a simultaneous sum-or-product game. The announced number is 35. Everyone says NO in the first round. What will happen in the second round?
\end{puzzle}

\begin{puzzle}[One person always knows.]\label{puzzle:6}
Three people are playing a simultaneous sum-or-product game. The announced number $M$ is semiprime. What is $M$, if someone always says YES in the first round independently of the distribution of the numbers?
\end{puzzle}

\begin{puzzle}[Half of us NO, and half of us don't.]\label{puzzle:10}
Four people Alice, Bob, Cindy, and Dylan have positive numbers on their foreheads. The emperor announces that the product or the sum of their numbers is 10. In the first round of the simultaneous game Alice and Bob said NO, and Cindy and Dylan said YES. What are the numbers of Alice and Bob?
\end{puzzle}

If Alice and Bob don't even know their numbers, how are we supposed to? This is the charm of common knowledge. 

We leave these puzzles for the readers and defer the solutions to Section~\ref{sec:solutions}.

\subsubsection{A complete analysis}

Let the announced number be $M=pq$, where $p$ and $q$ are distinct primes and $N > 2$.

As we discussed before there are several situations when one person, Alice, knows her number before anyone says anything:
\begin{itemize}
\item \textbf{It-must-be-the-product.} She knows $M$ is the product if the sum of the numbers she sees is $M$ or more. In particular, this happens if she sees $M$, or if she sees $p$ and $N > pq-p-1$, or if she sees $q$ and $N > pq - q -1$, or if she sees all ones and $N > M$.
\item \textbf{It-must-be-the-sum.} She knows $M$ is the sum if the product of the numbers she sees is not a factor of $M$. In particular, this happens if she sees a non-divisor of $M$, at least two $p$'s, or at least two $q$'s.
\item \textbf{Whatever.} She doesn't know whether $M$ is the sum of the product, but knows her number because both possibilities give the same answer. Suppose she sees one of the prime factors of $M$ not counting ones and $N = pq-p-q+2$, then she knows that her number must be the other prime factor. If she sees $p$ and $q$ and some ones, and again $N = pq-p-q+2$ she knows that her number must be 1.
\end{itemize}

We divide our discussion into three cases describing the set of numbers on foreheads: at least two non-divisors of $M$, one non-divisor of $M$, and no non-divisors of $M$.

\textbf{Case 1. At least two non-divisor of $M$.} If at least 2 people do not have divisors of $M$, then everyone sees at least one non-divisor. Everyone says YES in the first round, and the game is over. 

\textbf{Case 2. One non-divisor of $M$.} If exactly one person, Alice, does not have a proper divisor, then $N-1$ people say YES in the first round. In addition, if Alice sees more than one $p$ or more than one $q$, then Alice says YES and the game is over.

If Alice sees ones and not more than one $p$ and not more than one $q$, then she doesn't know whether $M$ is the product. As she has a proper divisor, there are two possibilities for what she can have. She says NO in the first round. We discuss subcases depending on what Alice sees.

\textbf{Subcase 2a.} Alice sees only ones. Then she doesn't know whether her number is $M$ or $M - N+1$, and she will never know her number.

\textbf{Subcase 2b.} Alice sees ones and one proper divisor of $M$, say $p$. Then her number is either $q$ or $M-p -N+2$. If her number was $q$, then the person with $p$ would have said NO in the first round. Thus Alice can deduce that $M$ is the sum and her number is not $q$. She says YES in the second round.

\textbf{Subcase 2c.} Alice sees ones and both $p$ and $q$. Alice figures out that if she had a one, then everyone else would have said NO in the first round. So she knows her number in the second round.

\textbf{Case 3. No non-divisors of $M$.} The last case is if everyone's number is in the set $\{1,p,q\}$. If there are 3 or more $p$'s or 3 or more $q$'s, everyone sees at least 2 $p$'s or 2 $q$'s, so everyone knows it's the sum in the first round and says YES. If there are 2 $p$'s and 2 $q$'s, each person still sees at least 2 $p$'s or 2 $q$'s, so, just as before, they all know that $M$ is the sum. Everyone says YES. We are left with the cases when either all numbers are ones or the non-ones form the following sets: $(p,p,q)$, $(p,q,q)$, $(p,q)$, $(p,p)$, $(q,q)$, $(p)$, $(q)$. We consider subcases.

\textbf{Subcase 3a.} The numbers are $p,p,q$ not counting ones (the case of $p,q,q$ is similar). Then $M=2p+q+N-3$ and people with $p$ do not know whether $M$ is the sum or the product. Everyone except the two people with $p$'s says YES in the first round. Then in round 2, since $pq \neq p+q+N-2$, each person with $p$ can say YES since if they had $1$, the person with $q$ would not have said YES in round 1. Thus two people with $p$ are the only people to say NO in the first round and they say YES in the second round.

\textbf{Subcase 3b.} The numbers are $p,p$ not counting ones (the case of $q,q$ is similar). Then $M=2p+N-2$ and people with $p$ do not know whether $M$ is the sum or the product. Everyone except the two people with $p$'s says YES in round 1. Each person with $p$ knows that if they had $q$, the people with $1$ would not have said YES in the first round since they didn't know whether they had $1$ or $M-p-q-(N-3)=p-q+1\neq 1$ since $p\neq q$. So both people with $p$ say YES in round 2, and the game is over.

\textbf{Subcase 3c.} The numbers are $p,q$ not counting ones. If $N \geq pq-p-q+2$, then everyone says YES in round 1. Otherwise, everyone says NO. 

Consider the second round when everyone says NO. It becomes common knowledge that all the numbers are proper divisors. People with ones know that they can't have $p$ as otherwise the person with $q$ would have said YES. Similarly, they can't have $q$. So they figure out that they have ones for the second round. 

Suppose Alice is the person with $p$. She knows that she can't have a non-divisor as, otherwise, everyone would have said YES in the first round. She reasons that she can't have $q$, as otherwise the people with ones would have said YES. She also can't have 1, since  then we would get $pq=N+q-1 < N+p+q-2$, which contradicts the results in the first round. Alice concludes that she must have $p$. With similar reasoning for the person with $q$, everyone say YES in the second round. 

\textbf{Subcase 3d.} The numbers are $p$ and ones (the case of $q$ is similar). Suppose Alice has $p$. It follows that there are $pq-p$ people with one. In the first round, Alice doesn't know whether she has $p$ or $pq$, and everyone else doesn't know whether they have 1 or $q$. Thus, everyone says NO in the first round. In the second round, Alice knows that if she had $pq$, then all other people would know their numbers in the first round. So she must have $p$ and says YES in the second round. 

Consider a particular person with a one, say Bob. He knows that he has 1 or $q$. Bob knows that if he had $q$, then the sum of the numbers that any other person with one sees would be $pq+q-2$. This is greater than or equal to $pq$, so if Bob had $q$, all the other people with ones would know that $M$ is the product and know their numbers in the first round. Since they didn't know that, Bob knows he must have 1 and says that in the second round, as do all other people with one. This means that everyone figures out their numbers in the second round.

\textbf{Subcase 3e.} All the numbers are ones. Everyone doesn't know whether they have 1 or $pq$. They say NO in the first round. If a person had $pq$, then everyone else would have said YES in the first round. Thus, they all say YES in the second round. 

To summarize, the following things could happen: Everyone says YES in the first round; One person says NO in the first round and YES in the second round; One person says NO in the first round and this person never figures out the number; Two people say NO in the first round and say YES in the second round. Everyone says NO in the first round and YES in the second round.

\section{Solutions to Puzzles Left for Readers}\label{sec:solutions}

\textbf{Solution to Puzzle~\ref{puzzle:ns-sim}. Near-sighted people.}
If $N \leq 3$, then all the people are neighbors of each other, so this is the original puzzle in Section~\ref{sec:intro}. The 

Suppose $N > 3$. Then the two people next to the person with the red hat will know their hat color in the first round. Then everyone knows that the person in between them has the red hat. From here, everyone can figure out their hat color unless there are two people that are both in between the two people who said YES in the first round, and this occurs if and only if $N=4$. In this scenario, the two people who said NO in the first round will not be able to deduce which one of them has a red hat. The answer to this puzzle is $N=4$.

\textbf{Solution to Puzzle~\ref{puzzle:9fs}. Nine far-sighted people.} We number the people from 1 to 9 in order starting with anyone. If any person sees zero or three red hats in the first round, they will say YES.
Therefore, one out of every three consecutive people has to have a red hat by
the Pigeonhole Principle. As a result, either Person 1, Person 4, and Person
7 have red hats, Person 2, Person 5, and Person 8 have red hats, or Person 3,
Person 6, and Person 9 have red hats. Everybody sees exactly two red hats,
so they know the exact configuration of the hats, so they will all say YES in the second round.

\textbf{Solution to Puzzle~\ref{puzzle:10ns}. Ten near-sighted people in a line.} Suppose everyone has 1 and Josh has 23. Then Aidan says YES. As this is the only case when Aidan can say YES, everyone would know their number after Aidan's announcement. They need 10 announcements: all of them YESes. Any other configuration will take more time. Also, you can't do better than everyone knowing their number in the first turn. Thus, this case requires the smallest number of turns.

\textbf{Solution to Puzzle~\ref{puzzle:2cons}. Two consecutive numbers.} At the beginning, Alice will know her number if and only if Bob has 0. On the second turn Bob knows that his number is positive. He can guess his number only if Alice's number is 0 or 1. Before the third turn Alice knows that her number is more than 1. She can only guess her number if Bob's number is 1 or 2. Before the fourth turn Bob knows that his number is more than 2. He can only guess his number if Alice's number is 2 or 3, and, correspondingly Bob's number is either 3 or 4. Given that Bob guessed his number and it is odd, we know that Bob's number is 3 and Alice's is 2.

\textbf{Solution to Puzzle~\ref{puzzle:mdsummary}. One, Zero, Two, Three.}
People with the maximum all have the same information, so they act the same way. Similarly, people with the minimum give the same answers, and also normal people give the same answer. As there are three groups that said YES, we know that the total number of people is six.

In the first round, one person will say YES if and only if $m=0$ and there is exactly one person with the maximum. After that it becomes common knowledge that there is a zero and the game becomes equivalent to the game with hats.
If the number of people with 0 is $r$, then these $r$ people will say YES during round $r+1$. Hence, we have two people with zeros. After the latest round everyone else knows they are normal. They only way they can figure out their numbers is if the maximum is 2. So the numbers must be 0, 0, 1, 1, 1, 2.

\textbf{Solution to Puzzle~\ref{puzzle:35}. Everyone Says NO.} The only distribution of numbers that would lead to no one knowing their number in the first round is $p$, $q$, and 1. Therefore, since no one says YES in the first round, the distribution can be determined. Thus, the answer to Puzzle~\ref{puzzle:35} is that everyone will say YES in the second round. By the way, their numbers are 1, 5, and 7.

\textbf{Solution to Puzzle~\ref{puzzle:6}. One person always knows.} We know that if numbers are $p,q,1$ everyone says NO in the first round unless $N = pq-p-q+2$. For $N = 3$ this equation is equivalent to $(p-1)(q-1) = 2$. It follows that $A = pq = 6$. Thus, the answer to Puzzle~\ref{puzzle:6} is 6.

\textbf{Solution to Puzzle~\ref{puzzle:10}.} By our analysis of the general case, two people say NO in the first round if the numbers are $(p,p,q,1)$, $(p,q,q,1)$, $(p,p,1,1)$, and $(q,q,1,1)$. If the sum or product is 10, the only possibility is $(2,2,5,1)$. Thus Alice and Bob have 2.

\textbf{Solution to Puzzle~\ref{puzzle:50circ}.} If Bob's number is not a divisor of 50, Alice knows that 50 must be the sum and, thus, she knows her number. If his number is 50, she knows that 50 must be the product, and thus, she knows her number. She does not know her number if and only if Bob's number is a proper divisor of 50, that is one of 1, 2, 5, 10, or 25. 

Similarly if Alice's number is not a proper divisor of 50, Bob would know his number. In addition, Bob knows that he has a proper divisor after Alice's first turn. The only way for 50 to be a sum of proper divisors is for both people to have 25. This means that if Alice's number is not 25, Bob knows that 50 is not the product. Thus, Alice must have 25. In this case, Bob's number is either 2 or 25 and he does not know his number.

\textbf{Solution to Puzzle~\ref{puzzle:50sim}. Half of us NO, and half of us don't.} Both have 25.

\section{Acknowledgements}

This project was done as part of PRIMES STEP, a program at MIT for middle school and 9th grade, which allows
students to try research. Tanya Khovanova is the mentor of the program. All other authors were in grades 7 through 9. We are grateful to PRIMES STEP for this opportunity.

\end{document}